\documentclass[reqno,preprint,12pt]{elsarticle}
\usepackage[active]{srcltx}
\usepackage{bbm}
\usepackage{tipa}
\usepackage{amsthm}
\usepackage{amssymb}
\usepackage{stackrel}
\usepackage{mathrsfs}
\usepackage{amsmath,amscd}
\usepackage{lscape}
\usepackage{longtable}
\usepackage{booktabs}
\usepackage{rotating}
\usepackage[active]{srcltx}
\usepackage{amsfonts}
\usepackage[centerlast]{caption}
\usepackage{latexsym}
\usepackage{amsmath,amsthm}
\usepackage{amssymb,amsmath}
\usepackage[all]{xy}
\usepackage{graphicx}
\usepackage[dvips]{geometry}
\newtheorem{dl}{Theorem}[section]
\newtheorem{tl}[dl]{Corollary}
\newtheorem{yl}[dl]{Lemma}
\newtheorem{dy}[dl]{Definition}

\newtheorem{xinzhi}[dl]{Proposition}

\numberwithin{equation}{section}

\newproof{pot1}{Proof of Theorem \ref{ma1}}
\newproof{pot2}{Proof of Theorem \ref{ma2}}
\def\qed{\hfill \rule{4pt}{7pt}}
\def\pf{\noindent {\sl Proof.~}}
\newcommand{\poq}[2]{(#1;q)_{#2}}

\usepackage[colorlinks=true,linkcolor=green,filecolor=brown,citecolor=green]{hyperref}
\usepackage{hyperref}
\usepackage{epsf}

\begin{document}
\title{Dual forms of the orthogonality relations of some classical $q$-orthogonal polynomials}
\author{Qi Chen$^{a,}$\fnref{fn0}}
\fntext[fn0]{E-mail address: chenqi2022@zjnu.edu.cn}
\address[P. R. China]{School of Mathematical Sciences, Zhejiang Normal  University,Jinhua 321004,~P.~R.~China}
\author{Jin Wang$^{b,}$\fnref{fn3,fn44}}
\fntext[fn3]{Supported by the Natural Science Foundation of Zhejiang Province under (Grant~No. LY24A010012) and the National Natural Science Foundation of China (Grant~No. 12001492)}
\fntext[fn44]{Corresponding author. E-mail address: jinwang@zjnu.edu.cn}
\address[P.R.China]{School of Mathematical Sciences, Zhejiang Normal  University,Jinhua 321004,~P.~R.~China}
\author{Xinrong Ma$^{c,}$\fnref{fn1,fn2}}
\fntext[fn1]{Supported by the National Natural Science Foundation of China (Grant~No.  11971341 and 12471315)}
\fntext[fn2]{E-mail address: xrma@suda.edu.cn.}
\address[P.R.China]{Department of Mathematics, Soochow University, Suzhou 215006, P.R.China}

\begin{abstract}
In this paper, by introducing new matrix operations and using a specific inverse relation, we  establish  the dual forms of the orthogonality relations
for some well-known discrete and continuous   $q$-orthogonal polynomials from the Askey-scheme such as the little and big $q$-Jacobi, $q$-Racah,  (generalized) $q$-Laguerre, as well as the Askey-Wilson polynomials. As one of the most interesting results, we show that the Askey-Wilson $q$-beta integral represented in terms of  the VWP-balanced $\,_8\phi_7$ series is  just a dual form of the orthogonality relation of the Askey-Wilson polynomials.
\end{abstract}

\maketitle
\thispagestyle{plain}
\section{Introduction}\label{sec1}

 In the past decades   various methods  of finding and proving summation and transformation formulas of
   hypergeometric series have been proposed.  One of
typical  methods is the inverse relations investigated
systematically  by J. Riordan \cite{riordan} at first time and
developed thoroughly by Ch. Krattenthaler, L. C. Hsu and W. C. Chu, as well as others. The reader might consult
\cite{4,car, chu1, 111, gouldhsu,kratt, milne, zengguo} and related references
therein.  Of these, it is particularly noteworthy that a series of research works \cite{chu1} by W. C. Chu  displays that  the inverse relations, as the celebrated WZ-method \cite{aeqb} does, is one of powerful tools to  hypergeometric series. Indeed, in the  rich world of summation and transformation formulas,  many results seem totally different in
form but are equivalent to each others up to inverse relations. 

As we will see later,  similar equivalences  occur often
 in  the orthogonality relations of orthogonal polynomials. Nevertheless, it has been lacking of systematic and  full investigation
so far. It is the problem that the present paper will attack. For a more
general development of the   theory of orthogonal polynomials, the reader might consult    \cite{ismailbook,ismailadd1,ismailadd2} by M. E. H. Ismail.  Regarding
applications of $q$-series to orthogonal polynomials, the reader is referred to \cite[Chapter 7]{10}     for  a good survey.

Our aim of this paper is
to associate the orthogonality
relations of discrete and continuous   $q$-orthogonal polynomials listed in the $q$-Askey-scheme \cite{schem} with some summation and transformation formulas of
basic hypergeometric series  via the use of inverse relations, regardless of the cost of tedious
calculations involved. As rewarding, the latter would at least, even
if it is not always new  to us, offer a possibly novel insight into
 discrete and continuous $q$-orthogonal polynomials.

\subsection{Some preliminaries on inverse relations}

At this stage, we had better give some explanations  on the method of inverse relations.  Generally speaking, the inverse relations is
always related to a pair of reciprocal (or inverse) relations. Recall that
     a inverse relation or a (matrix) inversion
formula in the context of combinatorial analysis,  following up \cite{kra},  is usually defined to be a pair of
{\sl infinite-dimensional
    lower-triangular} (in short, ILT) matrices $F=(f_{n,k})_{n,k \in \mathbb{N}}$ and  $G=(g_{n,k})_{n,k \in \mathbb{N}}$, denoted by $(F,G)$,
    over the complex field $\mathbb{C}$  such that
  \begin{align}
 \sum_{n\geq i\geq k}f_{n,i}g_{i,k}=\sum_{n\geq i\geq k}g_{n,i}f_{i,k}=\delta_{n,k}\,\,\mbox{for all}\ n,k\in \mathbb{N},\label{inversedef}
\end{align}
where $\delta$ denotes the usual Kronecker delta, $\mathbb{N}$ is
the set of nonnegative integers. Basic application of such a pair of inverse matrices
is that it provides a standard technique for deriving new summation
formulas from known ones. More precisely, assume that $(f_{n,k})$
and $(g_{n,k})$ are inverses of each other, then of course the
 transformation of two sequences $\{a_n\}_{n\geq 0}$ and $\{b_n\}_{n\geq 0}$
\begin{equation}\sum _{k=0} ^{n}f_{n,k}a_{k}=b_{n}\quad \mbox{is equivalent to}\quad \sum _{k=0}^{n}
g_{n,k}b_{k}=a_{n}.\label{4.1}\end{equation} In other words, once a
relation in (\ref{4.1}) is known, then another comes true as its \emph{dual form}. That is the routine way how to use the inverse relations.

Needless to say,  finding the inverse $(g_{n,k})$, often denoted by $(f_{n,k})^{-1}$, enables us to establish even more summation formulas,  which is the one of fundamental roles that the inverse relations play in the theory of hypergeometric series. This problem, as pointed
out by Ir. Gessel and D. Stanton
 \cite[p.175, \S 2]{111},  is equivalent to the Lagrange inversion
formula. As of today,  there has been a more general inverse relation
called the $(f,g)$-inversion formula \cite{0020}, since it
contains all previously known results such as the Gould-Hsu inversion formula
\cite{gouldhsu}, the Krattenthaler inversion formula \cite{kra}, as well as the
Warnaar inversion formula \cite{1000} as special cases.
\begin{dl}\mbox{\rm \cite[Theorem 1.3]{0020} }\label{math4} Let $f(x,y)$ and $g(x,y)$ be two arbitrary functions over ${\Bbb C}$ in variables
 $x,y$ and $\{b_i\}_{i\in \mathbb{N}}$ be such a
sequence that none of the terms $g(b_i,b_j)$  and $f(x_i,b_j)$ in
the denominator of the right  side of (\ref{a12}) and
(\ref{b12}) vanish. Suppose that $g(x,y)$ is antisymmetric, i.e.,
$g(x,y)=-g(y,x)$.
 Let $F=(f_{n,k})_{n,k\in \mathbb{N}}$ and $G=(g_{n,k})_{n,k\in \mathbb{N}}$ be two  ILT matrices
with the entries given by
\begin{subequations}\begin{align}
f_{n,k}&=\frac{\prod_{i=k}^{n-1}f(x_i,b_k)}
{\prod_{i=k+1}^{n}g(b_i,b_k)}\label{a12}\qquad\mbox{and}\\
g_{n,k}&=
\frac{f(x_k,b_k)}{f(x_n,b_n)}\frac{\prod_{i=k+1}^{n}f(x_i,b_n)}
{\prod_{i=k}^{n-1}g(b_i,b_n)},\label{b12}\quad\mbox{respectively}.
\end{align}\end{subequations}  Then $F=(f_{n,k})_{n,k\in \mathbb{N}}$ and $G=(g_{n,k})_{n,k\in
\mathbb{N}}$ is a pair of inverse relations if and only if  $ f(x,y)$ and $g(x,y)$ satisfy that for any $a,b,c,x\in
\mathbb{C}$,
\begin{align}
g(b,c)f(x,a)+g(c,a)f(x,b)+g(a,b)f(x,c)=0.\label{add}
\end{align}
\end{dl}
Up to now, there have been  found many  pairs of  $f(x,y)$ and $g(x,y)$  satisfying
(\ref{add}),  each  pair of which in turn yields  an $(f,g)$-inversion formula.

\subsection{A special $(f,g)$-inversion formula}
To facilitate  our  discussions, we need an equivalent
 statement of Theorem \ref{math4}. The reader is referred to \cite{wangjinthesis} for its detailed proof.
\begin{yl}\label{dl2}
 Let $\{x_i\}_{i\in \mathbb{N}}$ and $\{b_i\}_{i\in \mathbb{N}}$
 be arbitrary sequences over $\mathbb{C}$ such that $b_i's$ are pairwise
 distinct, $f(x,y)$ and $g(x,y)$ are subject to  \eqref{add}. Then \begin{align}
F_n=\sum_{k=0}^{n}G_k f(x_k,b_k)
\frac{\prod_{i=0}^{k-1}g(b_i,b_n)}{\prod_{i=1}^{k}f(x_i,b_n)}\label{27}
\end{align}
if and only if \begin{align}
G_n=\sum_{k=0}^{n}F_k\frac{\prod_{i=1}^{n-1} f(x_i,b_k)}
 {\prod_{i=0,i\neq k}^{n}g(b_i,b_k)}.\label{28}
\end{align}
\end{yl}

  This lemma covers  a very useful inverse relation on which our argument relies heavily.

\begin{xinzhi}\label{agar}
Let $\textsf{N}(a)$ be the ILT matrix with the $(n,k)$-th entry given by
\begin{subequations}
\begin{align}
[\textsf{N}(a)]_{n,k}=\frac{\poq{q^{-n},aq^{n}}{k}}{\poq{q,aq}{k}}q^k.
\label{defmatrix-II}\end{align}
 Then the $(n,k)$-th entry of its inverse $\textsf{N}^{-1}(a)$ must be
 \begin{align}
[\textsf{N}^{-1}(a)]_{n,k}=\frac{\poq{a,q^{-n}}{k}}{\poq{q,aq^{1+n}}{k}}\frac{1-aq^{2k}}{1-a}q^{kn}.
 \end{align}
 \end{subequations}
\end{xinzhi}
\pf It suffices to take
$f(x,y)=x-y,g(x,y)=x-y,b_i=q^{-i},x_i=aq^{i}$ in Lemma \ref{dl2}. Note that $g(x,y)=x-y $  and $f(x,y)=x-y$ satisfy \eqref{add}
 .
\qed

We remark that  Proposition  \ref{agar} is originally due to
 Carlitz \cite{car} and appeared in a different version
in the study of Bailey chains due to G. E. Andrews \cite[Chapter 3]{andrews}, which has now been known as a
$q$-analogue  of the famous Gould-Hsu inversion
formula\,\cite{gouldhsu}.

\subsection{Sketch of main idea and notation}

As mentioned earlier,  we will apply  Proposition \ref{agar}
  to find any possibly new  form of the orthogonality relation of
  discrete   or continuous $q$-orthogonal polynomials arising mainly from the Askey-scheme \cite{schem}. Our motivation is built on a common feature:  the
orthogonality relations of  many $q$-orthogonal polynomials listed in  the Askey-scheme  \cite{schem} are always  identified with  certain summation or transformation formulas. This universal phenomenon would not be  transparent, as far as we are aware, unless   viewed from the inverse relation \eqref{inversedef}.

Briefly, our argument can be  summarized  as follows:  first to write the
orthogonality relation of the given discrete orthogonal polynomials  as the
identity of the $(n,m)$-entries of certain matrix equation
\begin{align}
[A\left(v_i\delta_{i,j}\right)A^T]_{n,m}=[\left(h_i\delta_{i,j}\right)]_{n,m}.
\label{equI}
\end{align}
It is more often than not that the matrix $A$  can be expressed as
\begin{align} A=CX,\nonumber
\end{align}
while $C$ is of the form $D_1YD_2$,  $Y\in\{\textsf{N}(a), \textsf{N}^{-1}(a)\}$, $\textsf{N}(a)$ is given by Proposition \ref{agar}, $D_i(i=1,2)$ are diagonal. Thanks to this inverse relation, we are able to transform  the  matrix equation
(\ref{equI}) by the inverse relations to
\begin{align}
X\left(v_i\delta_{i,j}\right)X^T=C^{-1}\left(h_i\delta_{i,j}\right)C^{-T},\nonumber
\end{align}
hereafter,  the superscript $T$ denotes matrix transpose as usual. Accordingly, we have the relation often taken in form of
hypergeometric series
\begin{align}[X\left(v_i\delta_{i,j}\right)X^T]_{n,m}=[C^{-1}\left(h_i\delta_{i,j}\right)C^{-T}]_{n,m}.\label{equII}\end{align}
 This is precisely what we want. Henceforth,  (\ref{equII}) is called  the \emph{dual form} of  (\ref{equI}), since they are equivalent to each other.

The present paper, as further study of  the
 $(f,g)$-inversion formula \cite{0020,0021}, is organized as follows.  In  Section 2 we will apply the inverse relation $(\textsf{N}(a), \textsf{N}^{-1}(a))$ to find some dual forms of the orthogonality
relations of some discrete  $q$-orthogonal polynomials, among include the little $q$-Jacobi, $q$-Racah, and $q$-Laguerre orthogonal polynomials.  Section 3 is devoted to dual forms of the orthogonality
relations of some continuous $q$-orthogonal polynomials, such as the big
$q$-Jacobi and the Askey-Wilson orthogonal polynomials. To that end, we need to introduce a kind of new matrix operations. In both cases, some  new   summation and transformation formulas are obtained as byproducts.

  {\sl Notations and conventions.}
Throughout this paper we will use the standard notation and terminology
for basic hypergeometric series as in G. Gasper and M. Rahman
\cite{10}. For
instance, given a (fixed) complex number $q$ with $0<|q|<1$, a complex number
$a$ and a nonnegative integer $n$, denote the rising shifted factorial by
 \begin{align}\nonumber
   &(a;q)_{\infty}:=\prod_{k=0}^{\infty}(1-aq^k),\quad (a;q)_n:=
\frac{\poq{a}{\infty}}{\poq{aq^n}{\infty}}, \nonumber\\
   &(a_1,a_2,\cdots,a_m;q)_n = (a_1;q)_n(a_2;q)_n\cdots
   (a_m;q)_n.\nonumber
    \end{align}
 The {\em basic hypergeometric series\/} with the base $q$  and the argument $z$ is
defined by
\begin{equation}\nonumber{}_{r}\phi _{s}\left[\begin{matrix}a_{1},a_{2},\dots ,a_{r}\\ b_{1},b_{2},\dots ,b_{s}\end{matrix}
; q, z\right]:=\sum _{n=0} ^{\infty }\frac{\poq {a_{1},a_{2},\cdots
,a_{r}}{n}}{\poq
{q,b_{1},b_{2},\cdots,b_{s}}{n}}t_n(1+s-r)z^{n},\end{equation}
where the nonnegative integers  $r\leq s+1$ and $
t_n(k)=(-1)^{nk}q^{kn(n-1)/2}.$ In particular, the
${}_{r+1}\phi _{r}$ series is said to be very-well-poised (in brief, VWP), if all parameters $a_i$ and $b_i$ are
subject to the conditions that 
\begin{equation}a_1q=a_2b_1=a_3b_2=\cdots=a_{r+1}b_r; a_2=q\sqrt{a_1},a_3=-q\sqrt{a_1}.\nonumber
\end{equation}
Such  ${}_{r+1}\phi _{r}$ series is often denoted by the more compact
notation
\begin{equation}{}_{r+1}W_{r}(a_1;a_4,\dots,a_{r+1};q, z).\nonumber\end{equation}
In addition, for any sequence $\left\{A_n\right\}_{n\in{\mathbb{Z}}}$, we will employ the convention of defining
 (cf. \cite[(3.6.12)]{10})
\begin{align}\nonumber
\prod_{i=k}^nA_i:= \left\{\begin{array}{ll}A_{k}A_{k+1}\cdots A_{n},&\mbox{if}\,\, n\geq k;\\
1,&\mbox{if}\,\, n=k-1;\\
 1/(A_{n+1}A_{n+2}\cdots A_{k-1}),&\mbox{if}\,\, n\leq k-2
\end{array}
\right.
\end{align}
over the set of integers ${\Bbb Z}$ , the notation $A=(a_{n,k})$  to denote the
infinite  matrix $A$ with the $(n,k)$-th entry
$a_{n,k}$, or in short, $[A]_{n,k}=a_{n,k}.$
\section{Dual forms of the orthogonality relations of  discrete  $q$-Jacobi, $q$-Racah, and  $q$-Laguerre polynomials}\label{sec2}
As planned, in this part  we will apply  Proposition \ref{agar}
  to pursuit for any possibly new equivalent forms of the orthogonality relations of
  discrete  $q$-orthogonal polynomials from the Askey-scheme \cite{schem}.
We content ourselves only with displaying a few typical discrete  $q$-orthogonal polynomials from \cite{schem} as good examples to illustrate our idea described in Section 1.3.

\subsection{The little $q$-Jacobi
polynomials} Recall that  the little $q$-Jacobi polynomials
\cite[(3.12.1)]{schem} are defined, for
$a\not\in\{q^m : m \in \mathbb{Z},m<0\} $, by
\begin{align}
p_n(x;a,b;q):=\,_2\phi_1\left[%
\begin{array}{cc}
  q^{-n}, & abq^{n+1}\\
   & aq\\
\end{array};q,xq
\right].\label{jacobi}
\end{align}
\begin{yl}{\rm (cf. \cite[(3.12.2)]{schem})} The little $q$-Jacobi polynomials satisfy the orthogonality relation
\begin{align}
\sum_{k=0}^\infty\frac{a^kq^k\poq{bq}{k}}{\poq{q}{k}}p_n(q^k;a,b;q)p_m(q^k;a,b;q)
=h_n\delta_{n,m},\label{orthog}
\end{align}
where
\begin{equation}
h_n=\frac{\poq{abq^{2}}{\infty}}{\poq{aq}{\infty}}
\frac{\poq{q,bq}{n}}{
\poq{aq,abq}{n}}\frac{(1-abq)(aq)^n}{1-abq^{2n+1}}.
\nonumber\end{equation}
\end{yl}
We refer the reader to \cite{10} for an alternate proof of
(\ref{orthog}) by appealing to the nonterminating $q$-Saalsch\"{u}tz
 summation formula  (\cite[(II.24)]{10}). From the viewpoint of inverse relations, it is equivalent to
\begin{dl}\label{dl11} The dual form of the orthogonality relation  (\ref{orthog}) is  the terminating  $\,_{6}\phi_{5}$ summation formula (\cite[(II.20)]{10})
\begin{align}
\,_{6}W_{5}\left(abq;bq,q^{-n},q^{-m};q,aq^{1+m+n}\right)=\frac{\poq{abq^2}{n}\poq{abq^2}{m}\poq{aq}{n+m}}{\poq{aq}{n}\poq{aq}{m}\poq{abq^2}{n+m}}.
\label{32}    \end{align}
    \end{dl}
\pf
 It suffices to express  (\ref{orthog}) in terms of matrices. For this, we
introduce three infinite  matrices $A,B,X$ with the entries
\begin{align}
[A]_{n,k}=p_n(q^k;a,b;q),
[B]_{n,k}=\frac{a^kq^k\poq{bq}{k}}{\poq{q}{k}}\delta_{n,k},
[X]_{n,k}=\frac{\poq{abq^2}{n}}{\poq{aq}{n}}q^{nk}.\nonumber
\end{align}
By doing so, we can recast the orthogonality relation
(\ref{orthog})  in terms of these matrices  as
 \begin{align}
 ABA^T=\left(h_n\delta_{n,m}\right).\nonumber
\end{align}
From  (\ref{jacobi}) it is clear that
\begin{align}
 A=\textsf{N}(abq)X.\label{defA}
 \end{align}
 A substitution of  \eqref{defA} yields
\begin{align}\textsf{N}(abq)\left\{XBX^T\right\}\textsf{N}(abq)^T=\left(h_n\delta_{n,m}\right).\nonumber
\end{align}
By using the inverse relation $(\textsf{N}(a), \textsf{N}^{-1}(a))$, we get
\begin{align}
XBX^T=\textsf{N}^{-1}(abq)\left(h_n\delta_{n,m}\right)(\textsf{N}^{-1}(abq))^T.\label{ttt}
\end{align}
Upon equating  the $(n,m)$-th entries on both sides of (\ref{ttt}), we obtain
\begin{align}
&\displaystyle \sum_{k\geq 0}
\frac{\poq{abq^2}{n}}{\poq{aq}{n}}q^{kn}\frac{a^kq^k\poq{bq}{k}}{\poq{q}{k}}\frac{\poq{abq^2}{m}}{\poq{aq}{m}}q^{km}\nonumber
&\\
&=\displaystyle \sum_{k\geq 0}
\frac{\poq{q^{-n}}{k}}{\poq{abq^{2+n}}{k}}\frac{\poq{abq}{k}}{\poq{q}{k}}
\frac{1-abq^{2k+1}}{1-abq}\,q^{kn}\,\frac{\poq{abq^{2}}{\infty}}{\poq{aq}{\infty}}
\frac{(1-abq)(aq)^k}{1-abq^{2k+1}}\nonumber&\\
&\times\displaystyle \frac{\poq{q,bq}{k}}{ \poq{aq,abq}{k}}\,
\frac{\poq{q^{-m}}{k}}{\poq{abq^{2+m}}{k}}\frac{\poq{abq}{k}}{\poq{q}{k}}\frac{1-abq^{2k+1}}{1-abq}q^{km}.\nonumber&
\end{align}
After  series rearrangement and simplification,  noting that the sum on the left-hand side of the above identity can be evaluated in closed form by the $q$-binomial theorem \cite[(II.3)]{10}, it reduces
to
\begin{align}
&\sum_{k\geq 0}
\frac{\poq{q^{-n},q^{-m},abq,bq}{k}}{\poq{abq^{2+n},abq^{2+m},q,aq}{k}}
\frac{(1-abq^{2k+1})\left(aq^{1+n+m}\right)^k}{1-abq}\nonumber\\
&\qquad\qquad=
\frac{\poq{abq^2}{n}\poq{abq^2}{m}\poq{aq}{n+m}}{\poq{aq}{n}\poq{aq}{m}\poq{abq^2}{n+m}}.\nonumber
\end{align}
It is consistent with \eqref{32} after recast in standard  notation of  basic hypergeometric series.\qed

\subsection{The $q$-Racah polynomials} Another kind of  orthogonal polynomials  deserving our attention is the
$q$-Racah polynomials. Recall that the
$q$-Racah polynomials \cite[(7.2.17)]{10} is defined,  for
$N\geq m,n$, by
\begin{align}
R_n(x;a,b,c,N;q):=\,_4\phi_3\left[%
\begin{array}{cccc}
  q^{-n}, & abq^{n+1},&q^{-x},&cq^{x-N}\\
   & aq,&q^{-N},&bcq\\
\end{array};q,q
\right].\label{racah}
\end{align}
\begin{yl}{\rm (cf. \cite[(7.2.18)-(7.2.20)]{10})} The $q$-Racah polynomials satisfy the orthogonality relation
\begin{align}
\sum_{k=0}^\infty\rho(k;q)R_n(k;a,b,c,N;q)R_m(k;a,b,c,N;q)=\frac{\delta_{n,m}}{h_n},\label{738orthog1}
\end{align}
where
\begin{align}
\rho(k;q)=\frac{\poq{cq^{-N}}{k}}{\poq{q}{k}}\frac{(1-cq^{2k-N})}{(1-cq^{-N})}
\frac{\poq{aq,bcq,q^{-N}}{k}}{\poq{ca^{-1}q^{-N},b^{-1}q^{-N},cq}{k}}(abq)^{-k}\nonumber
\end{align}
and
\begin{align}
h_n=\frac{\poq{bq,aq/c}{N}} {\poq{abq^2,1/c}{N}}
\frac{(1-abq^{2n+1})}{(1-abq)}\frac{\poq{abq,aq,bcq,q^{-N}}{n}}
{\poq{q,bq,aq/c,abq^{N+2}}{n}}\left(\frac{q^N}{c}\right)^n.\nonumber
\end{align}
\end{yl}
As is to be expected, this orthogonality relation leads us to a new transformation between two
terminating $_8\phi_7$ series, stated as below.
\begin{dl}\label{dl12} The dual form of the orthogonality relation  (\ref{738orthog1}) is
   \begin{align}
& \,_{8}W_{7}\left(abq;bq,
aq/c,abq^{N+2},q^{-n},q^{-m};q,cq^{m+n-N}\right)\nonumber
\\
&=\frac{\poq{cq^{n-N}/a,q^{n-N}/b}{N-n}\poq{abq^2,q^{-n},cq^{n-N}}{m}}
{\poq{q^{-1-N}/ab,cq^{2n-N+1}}{N-n}\poq{aq,q^{-N},bcq}{m}}\label{433}\\
&\times\,_{8}W_{7}\left(
 cq^{2n-N};  cq^{m+n-N},q^{1+n}, aq^{1+n},
 bcq^{1+n},q^{n-N};q,\frac{1}{abq^{1+n+m}}\right),\nonumber
    \end{align}
 where  $N\geq n\geq m>0$.
       \end{dl}
 \pf
  Observe first that
the orthogonality relation of the $q$-Racah polynomials, as already given by
(\ref{738orthog1}), can be  reformulated  in terms of matrix operations as
\begin{align}
ABA^{T}=(\delta_{n,m}/h_n),\label{former}
\end{align}
where, we  introduce three infinite matrices $A,B,X$ such that
\begin{align}
&[A]_{n,k}=R_n(k;a,b,c,N;q),\,\,
[B]_{n,k}=\rho(k;q)\delta_{n,k},\,\,
[X]_{n,k}=\frac{\poq{abq^2,q^{-k},cq^{k-N}}{n}}{\poq{aq,q^{-N},bcq}{n}}.\nonumber
\end{align}
 Using these matrices and the  ILT matrix $\textsf{N}(a)$,  we easily restate  (\ref{racah})
  by the matrix identity
\begin{equation}
 A=\textsf{N}(abq)X.
 \nonumber\end{equation}
In the sequel, we are able to reformulate the orthogonality relation (\ref{738orthog1}) as
\begin{align}
XBX^T=\textsf{N}^{-1}(abq)\left(\frac{\delta_{n,m}}{h_n}\right)\textsf{N}^{-T}(abq).
\label{ttt0}
\end{align}
Upon equating the $(n,m)$-th entries on both sides, we obtain
\begin{align}
&\displaystyle \displaystyle \sum_{k\geq 0}
\frac{\poq{abq^2,q^{-k},cq^{k-N}}{n}}{\poq{aq,q^{-N},bcq}{n}}
\frac{\poq{abq^2,q^{-k},cq^{k-N}}{m}}{\poq{aq,q^{-N},bcq}{m}}\nonumber&\\
&\displaystyle\times
\frac{\poq{cq^{-N},aq,bcq,q^{-N}}{k}}{\poq{q,ca^{-1}q^{-N},b^{-1}q^{-N},cq}{k}}
\frac{(1-cq^{2k-N})(abq)^{-k}}{1-cq^{-N}}\nonumber
&\\
&=\displaystyle \displaystyle \sum_{k\geq 0}
\frac{\poq{q^{-n}}{k}}{\poq{abq^{2+n}}{k}}\frac{\poq{abq}{k}}{\poq{q}{k}}\frac{1-abq^{2k+1}}{1-abq}q^{kn}
\frac{\poq{q^{-m}}{k}}{\poq{abq^{2+m}}{k}}\frac{\poq{abq}{k}}{\poq{q}{k}}\frac{1-abq^{2k+1}}{1-abq}q^{km}\nonumber&\\
&\displaystyle\times  \frac{\poq{abq^2,1/c}{N}}{\poq{bq,aq/c}{N}}
\frac
{\poq{q,bq,aq/c,abq^{N+2}}{k}}{\poq{abq,aq,bcq,q^{-N}}{k}}\frac{(1-abq)}{(1-abq^{2k+1})}\left(\frac{c}{q^N}\right)^k.\nonumber&
\end{align}
On considering    the factors $\poq{q^{-k}}{n}\poq{q^{-k}}{m}=0$ for $k<n$ or $k<m$, we now assume that $k\geq n\geq m$
without any loss of generality. In this case, change the  index
 of summation from $k$ to $k+n$ on the left-hand and reformulate the last identity  in the form
\begin{align}
&\frac{\poq{abq^2}{n}\poq{abq^2}{m}}{\poq{aq,q^{-N},bcq}{n}\poq{aq,q^{-N},bcq}{m}}
\frac{\poq{cq^{1-N},aq,bcq,q^{-N}}{n}}{\poq{q,ca^{-1}q^{-N},b^{-1}q^{-N},cq}{n}}
\frac{(abq)^{-n}}{(1-cq^{n-N})}\nonumber
\\
&\sum_{k\geq 0}
\poq{q^{-k-n},cq^{k+n-N}}{n}\poq{q^{-k-n},cq^{k+n-N}}{m}
(1-cq^{2k+2n-N})\nonumber\\
&\qquad\qquad\times\frac{\poq{cq^{n-N}}{k}}{\poq{q^{1+n}}{k}}\frac{\poq{aq^{1+n},bcq^{1+n},q^{n-N}}{k}}
{\poq{ca^{-1}q^{n-N},b^{-1}q^{n-N},cq^{1+n}}{k}}(abq)^{-k}\nonumber
\\
&=\frac{\poq{abq^2,1/c}{N}}{\poq{bq,aq/c}{N}}\sum_{k\geq 0}
\frac{\poq{q^{-n}}{k}}{\poq{abq^{2+n}}{k}}\frac{\poq{abq}{k}}{\poq{q}{k}}\frac{1-abq^{2k+1}}{1-abq}
\frac{\poq{q^{-m}}{k}}{\poq{abq^{2+m}}{k}}\frac{\poq{abq}{k}}{\poq{q}{k}}\nonumber\\
&\qquad\qquad\times
\frac{1-abq^{2k+1}}{1-abq}\frac{(1-abq)}{(1-abq^{2k+1})}\frac{\poq{q}{k}}{\poq{abq}{k}}\frac
{\poq{bq,aq/c,abq^{N+2}}{k}}{\poq{aq,bcq,q^{-N}}{k}}\left(cq^{m+n-N}\right)^k\nonumber
.
\end{align}
At this stage, in view of  the basic relation \cite[(I.13)]{10}
\begin{align}
\poq{xq^{-k}}{n}=\frac{\poq{x}{n}\poq{q/x}{k}}{\poq{q^{1-n}/x}{k}}
q^{-nk},\label{basic}
\end{align} it is easy to check that for integer $n\geq 0$, there hold
    \begin{align} &\poq{q^{-k-n}}{n}
=\frac{\poq{q^{-n}}{n}\poq{q^{1+n}}{k}}{\poq{q}{k}}q^{-nk}\quad\mbox{and}\nonumber\\
&\poq{cq^{k+n-N}}{n}=\frac{\poq{cq^{n-N}}{n}\poq{cq^{2n-N}}{k}}{\poq{cq^{n-N}}{k}}.\nonumber
\end{align}
Simplify the preceding identity by  these expressions. The result is
\begin{align}
&\frac{\poq{abq^2,q^{-n},cq^{n-N}}{m}\poq{abq^2}{n}}
{\poq{aq,q^{-N},bcq}{m}\poq{ca^{-1}q^{-N},b^{-1}q^{-N}}{n}}
\frac{\poq{cq^{1-N}}{2n}(-ab)^{-n}q^{-n(n+3)/2}}{\poq{cq}{n}}\nonumber\\
&\sum_{k\geq 0}
\frac{\poq{cq^{2n-N},cq^{m+n-N},q^{1+n},aq^{1+n},bcq^{1+n},q^{n-N}}{k}}{\poq{q,q^{1-m+n},cq^{n-N},ca^{-1}q^{n-N},
b^{-1}q^{n-N},cq^{1+n}}{k}}
\frac{(1-cq^{2k+2n-N})}{(1-cq^{2n-N})}(abq^{1+m+n})^{-k}\nonumber
\\
&=\frac{\poq{abq^2,1/c}{N}}{\poq{bq,aq/c}{N}}\sum_{k\geq 0} \frac
{\poq{abq,q^{-n},q^{-m},bq,aq/c,abq^{N+2}}{k}}{\poq{q,abq^{2+n},abq^{2+m},aq,bcq,q^{-N}}{k}}\frac{1-abq^{2k+1}}{1-abq}\left(cq^{m+n-N}\right)^k\nonumber
,
\end{align}
which,  written out in standard notation of basic hypergeometric series, is recognized to be
 \begin{align}
& \,_{8}W_{7}\left(abq;bq,
aq/c,abq^{N+2},q^{-n},q^{-m};q,cq^{m+n-N}\right)\label{eq244}
\\
&=C_0\,_{8}W_{7}\left(
 cq^{2n-N};  cq^{m+n-N},q^{1+n}, aq^{1+n}, bcq^{1+n},q^{n-N};q,\frac{1}{abq^{1+n+m}}\right), \nonumber
    \end{align}
  where
 \begin{align}C_0&=\frac{\poq{bq,aq/c}{N}}{\poq{cq^{-N}/a,q^{-N}/b}{n}}\frac{\poq{abq^2,q^{-n},cq^{n-N}}{m}}{\poq{aq,q^{-N},bcq}{m}}\nonumber
\\
&\times\frac{\poq{abq^2}{n}\poq{cq^{1-N}}{N}\poq{cq}{2n-N}}{\poq{abq^2,1/c}{N}\poq{cq}{n}}
\left(\frac{-q^{-(n+3)/2}}{ab}\right)^{n}.\nonumber
 \end{align}
Note that the $C_0$ can be further simplified  by (\ref{basic}) to
\begin{align}C_0=\frac{\poq{cq^{n-N}/a,q^{n-N}/b}{N-n}\poq{abq^2,q^{-n},cq^{n-N}}{m}}
{\poq{q^{-1-N}/ab,cq^{2n-N+1}}{N-n}\poq{aq,q^{-N},bcq}{m}}.\label{C0}
 \end{align} The desired result follows from \eqref{eq244} at once.\qed

 It is of interest to see that  by applying Watson's ($q$-Whipple) transformation  (cf. \cite[(III.17)]{10})
    to  two $_8\phi_7$ sums on both sides of (\ref{433}),  we recover a finite form of   Sear's
    $\,_{4}\phi_{3}$ transformation formula \cite[(III.16)]{10}.
\begin{tl} For integers $N\geq n\geq m>0$, there holds
\begin{align}
&\,_{4}\phi_{3}\left[%
\begin{array}{cccc}
q^{-m},& aq/c,& q^{-n},&q^{-1-N}/b\\
&q^{-N},& aq,&q^{-m-n}/bc
\end{array};q,q\right]\label{2.14}\\
&=\frac {\poq{abq^{2+n},q^{-n},cq^{n-N}}{m}}{
\poq{bcq^{1+n},q^{-N},aq} {m}} \,_{4}\phi_{3}\left[%
\begin{array}{cccc}
q^{-m},& aq^{1+n},& bcq^{1+n},&q^{n-N}\\
&q^{1+n-m},& cq^{n-N},&abq^{2+n}
\end{array};q,q\right].\nonumber
    \end{align}
\end{tl}\pf As indicated above, by applying Watson's transformation \cite[(III.17)]{10}  to  the $\,_{8}\phi_{7}$
 series on both sides of (\ref{433}) simultaneously, we obtain
\begin{align}
&\nonumber\frac{\poq{abq^2,bcq^{1+n}} {m}}
{\poq{bcq,abq^{2+n}}{m}}\,_{4}\phi_{3}\left[%
\begin{array}{cccc}
q^{-m},& aq/c,& q^{-n},&q^{-1-N}/b\\
&\nonumber q^{-N},& aq,&q^{-m-n}/bc
\end{array};q,q\right] \\
&\nonumber=C_0\frac{\poq{cq^{1+2n-N},q^{-1-N}/ab} {N-n}}
{\poq{cq^{n-N}/a,q^{n-N}/b}{N-n}}\,_{4}\phi_{3}\left[%
\begin{array}{cccc}
q^{-m},& aq^{1+n},& bcq^{1+n},&q^{n-N}\\
&q^{1+n-m},& cq^{n-N},&abq^{2+n}
\end{array};q,q\right]\nonumber,
    \end{align}
   where
$C_0
$ is given by \eqref{C0}.
To simplify further, we reformulate it as
 \begin{align}
&\,_{4}\phi_{3}\left[%
\begin{array}{cccc}
q^{-m},& aq/c,& q^{-n},&q^{-1-N}/b\\
&q^{-N},& aq,&q^{-m-n}/bc
\end{array};q,q\right]\label{000}\\
&=C_1\,_{4}\phi_{3}\left[%
\begin{array}{cccc}
q^{-m},& aq^{1+n},& bcq^{1+n},&q^{n-N}\\
&q^{1+n-m},& cq^{n-N},&abq^{2+n}
\end{array};q,q\right]\nonumber
    \end{align}
  by defining
    \begin{align}\nonumber C_1&=\frac
{\poq{bcq,abq^{2+n}}{m}}{\poq{abq^2,bcq^{1+n}}
{m}}\frac{\poq{cq^{1+2n-N},q^{-1-N}/ab} {N-n}}
{\poq{cq^{n-N}/a,q^{n-N}/b}{N-n}}\\
&\nonumber\times\frac{\poq{cq^{n-N}/a,q^{n-N}/b}{N-n}\poq{abq^2,q^{-n},cq^{n-N}}{m}}
{\poq{q^{-1-N}/ab,cq^{2n-N+1}}{N-n}\poq{aq,q^{-N},bcq}{m}}.
 \end{align}
A direct simplification yields
 \begin{align}\nonumber C_1=\frac
{\poq{abq^{2+n},q^{-n},cq^{n-N}}{m}}{ \poq{bcq^{1+n},q^{-N},aq}
{m}}.
 \end{align}
 Thus the transformation \eqref{2.14}  follows from (\ref{000}).\qed

  \subsection{The $q$-Laguerre polynomials}\label{laguerre}
A well-known version of the $q$-Laguerre polynomials
\cite[(3.21.1)]{schem} is given, for $\alpha>-1$, by
\begin{align}
L_n^{(\alpha)}(x;q):=\frac{1}{\poq{q}{n}}\,_2\phi_1\left[%
\begin{array}{cc}
  q^{-n}, & -x\\
   & 0\\
\end{array};q,q^{n+\alpha+1}
\right].\label{lag}
\end{align}
\begin{yl}{\rm (cf. \cite[(3.21.3)]{schem})} The $q$-Laguerre polynomials
 satisfy the following discrete orthogonality relation
\begin{align}
\sum_{k=-\infty}^\infty\frac{q^{k+k\alpha}}{\poq{-cq^k}{\infty}}L_n^{(\alpha)}(cq^k;q)
L_m^{(\alpha)}(cq^k;q)=h_n\delta_{n,m},\label{orthog3}
\end{align}
where $c>0$ and
\begin{align}
h_n=\frac{\poq{q,-cq^{\alpha+1},-c^{-1}q^{-\alpha}}{\infty}}{\poq{q^{\alpha+1},-c,-c^{-1}q}{\infty}}
\frac{\poq{q^{\alpha+1}}{n}}{ \poq{q}{n}}q^{-n}. \label{orthog3-333}
\end{align}
\end{yl}
Evidently, this orthogonality relation is distinct from  all preceding cases. Viewed in light of inverse relations, it implies
\begin{dl}\label{dl1133}For $\alpha>-1,
c>0$, the dual form  of the orthogonality relation (\ref{orthog3}) is    the
 basic identity
 \begin{align}
\sum_{k=0}^n\begin{bmatrix}n\\k\end{bmatrix}_qy^{n-k}x^k\poq{y}{k}=
\sum_{k=0}^n\frac{\poq{q^{-n},x,y}{k}}{\poq{q}{k}}q^k,
\end{align}
where  $(x,y)=(q^{-m},q^{\alpha+1})$, the $q$-binomial coefficient
\begin{equation}\nonumber\begin{bmatrix}n\\k\end{bmatrix}_q:=\frac{\poq{q}{n}}{\poq{q}{k}\poq{q}{n-k}}.
\end{equation}
    \end{dl}
\pf Proceeding as before, in order to express  (\ref{orthog3}) in terms of matrix operations, now we
define  three infinite  matrices $A,B,X$ with the $(n,k)$-th entries
\begin{align}\nonumber
[A]_{n,k}=L_n^{(\alpha)}(cq^k;q),
[B]_{n,k}=\frac{q^{k+k\alpha}}{\poq{-cq^k}{\infty}}\delta_{n,k},
[X]_{n,k}=\poq{-cq^{k}}{n}q^{(\alpha+1)n}.
\end{align}
Note that in this case,  the index $k\in {\Bbb Z}$.
From  (\ref{lag}) it is clear that
\begin{align}\nonumber
 A=\left(\frac{\delta_{n,k}}{\poq{q}{n}}\right)\textsf{N}^{-1}(0)X.
 \end{align}
Now, with these facts and  $\textsf{N}(a)$
 given by (\ref{defmatrix-II}), the orthogonality relation (\ref{orthog3}) of the $q$-Laguerre
polynomials can be reformulated as
\begin{align}
XBX^T=\textsf{N}(0)\left(\poq{q}{n}\delta_{n,k}\right)\left(h_n\delta_{n,m}\right)\left(\poq{q}{n}\delta_{n,k}\right)\textsf{N}(0)^{T}.\label{ttts}
\end{align}
By equating the $(n,m)$-th entries on both sides of (\ref{ttts}) gives
\begin{align}\nonumber
\displaystyle \sum_{k=-\infty}^\infty
\poq{-cq^{k}}{n}q^{(\alpha+1)n}\,\frac{q^{k+k\alpha}}{\poq{-cq^k}{\infty}}\,\poq{-cq^{k}}{m}q^{(\alpha+1)m}
=\displaystyle \sum_{k=0}^\infty \poq{q^{-m}}{k}\poq{q^{-n}}{k}h_kq^{2k}.&
\end{align}
Upon substituting \eqref{orthog3-333} for $h_n$ and simplifying, we obtain
\begin{align}
&\displaystyle q^{(\alpha+1)(n+m)}\sum_{k=-\infty}^\infty
\frac{\poq{-cq^{k}}{m}\poq{-cq^{k}}{n}}{\poq{-cq^{k}}{\infty}}q^{k(1+\alpha)}\label{wangjinwenti}
\\
&=\displaystyle\frac{\poq{q,-cq^{\alpha+1},-c^{-1}q^{-\alpha}}{\infty}}{\poq{q^{\alpha+1},-c,-q/c}{\infty}}
\sum_{k=0}^{\min\{m,n\}} \frac{\poq{q^{-n},q^{-m},q^{\alpha+1}}{k}}{
\poq{q}{k}}q^{k}.\nonumber
\end{align}
Next, by applying the
$q$-binomial theorem \cite[(II.3)]{10} to $\poq{-cq^k}{n}$  on the left-hand side  of \eqref{wangjinwenti} and  interchanging the order of summation,  we obtain
\begin{align}\mbox{LHS of  \eqref{wangjinwenti}}
&\nonumber=\displaystyle~q^{(\alpha+1)(n+m)}\frac{\poq{-c}{m}}{\poq{-c}{\infty}}\sum_{i=0}^n\begin{bmatrix}n\\i\end{bmatrix}_qq^{\binom{i}{2}}c^i\sum_{k=-\infty}^\infty
\poq{-cq^{m}}{k}q^{k(1+i+\alpha)}\\
&\nonumber=q^{(\alpha+1)(n+m)}\frac{\poq{-c}{m}}{\poq{-c}{\infty}}\sum_{i=0}^n\begin{bmatrix}n\\i\end{bmatrix}_qq^{\binom{i}{2}}c^i
\frac{\poq{q,-cq^{m+1+i+\alpha},-q^{-m-i-\alpha}/c}{\infty}}{\poq{-q^{1-m}/c,q^{1+i+\alpha}}{\infty}}.
\end{align}
The last equality comes from Ramanujan's $\,_1\psi_1$ summation
formula \cite[(II.29)]{10} with $b=0$.  Some routine  simplifications on this expression   finally reduce \eqref{wangjinwenti} to the special case $x\to
q^{-m},y\to q^{\alpha+1}$ of the
 basic identity
 \begin{align}
\sum_{k=0}^n\begin{bmatrix}n\\k\end{bmatrix}_qy^{n-k}x^k\poq{y}{k}=
\sum_{k=0}^n\frac{\poq{q^{-n},x,y}{k}}{\poq{q}{k}}q^k.\nonumber 
\end{align}
It is just a limiting case of \cite[(III.13)]{10}. Hence the theorem is proved.\qed

It  may come as a surprise that  by the above  argument, we can extend  the orthogonality relation \eqref{orthog3} to the following
\begin{dl} For $|y|<1$,  we redefine the $q$-Laguerre polynomials by \begin{align}
L_n(x,y;q):=\frac{1}{\poq{q}{n}}\,_2\phi_1\left[%
\begin{array}{cc}
  q^{-n}, & -x\\
   & 0\\
\end{array};q,yq^{n}
\right].\label{lagg}
\end{align}
Then,  for $c>0$, $L_n(x,y;q)$  satisfies the following orthogonality relation:
\begin{align}
\sum_{k=-\infty}^\infty\frac{y^{k}}{\poq{-cq^k}{\infty}}L_n(cq^k,y;q)L_m(cq^k,y;q)=\frac{\poq{q,-cy,-q/cy}{\infty}}{\poq{y,-c,-c^{-1}q}{\infty}}
\frac{\poq{y}{n}}{ \poq{q}{n}}q^{-n}\delta_{n,m}.\label{laggorth}
\end{align}
\end{dl}
Evidently, $L_n(x,y;q)$ is  polynomial in two variables $x$ and
$y$. Thus  we may guess that it should be an orthogonal  polynomial in
 $y$ with respect to  certain measure, too.
 \section{Dual forms of the orthogonality relations of   the continuous big $q$-Jacobi and Askey-Wilson polynomials}\label{sec3}
As a matter of fact,  the foregoing argument  to find dual forms of the orthogonality relations  by  the inverse relations is also valid for continuous $q$-orthogonal polynomials. For this purpose, we need a few  new concepts  which are analogue to the usual matrix algebras or the incidence algebras over partially ordered set \cite{stanleybook}.
 \subsection{New matrix operations}
 \begin{dy} \label{newmatrix0} Given two finite or infinite index sets $\mathbb{K}_1, \mathbb{K}_2 $ and the number field $\mathbb{K}$,
 any binary mapping
\begin{equation}\nonumber
 f: \mathbb{K}_1\times \mathbb{K}_2\to  \mathbb{K}
 \end{equation}
 is called a two-dimensional matrix over $\mathbb{K}_1\times\mathbb{K}_2 $. The set of such matrices is denoted by
\begin{equation}\nonumber\Omega^{\mathbb{K}_1\times \mathbb{K}_2}:=\left\{A=(f(x,y))\big| x\in\mathbb{K}_1, y\in \mathbb{K}_2 \right\}.\end{equation}
\end{dy}
As customary,   we employ  the notation $[A]_{x,y}$ for its entry $f(x,y)$.
\begin{dy} \label{newmatrix}
 For arbitrary matrices  $X\in \Omega^{\mathbb{Z}\times \mathbb{C}}$, $Y\in \Omega^{\mathbb{C}\times \mathbb{Z}}$, and $A\in \Omega^{\mathbb{Z}\times \mathbb{Z}}$, we  further introduce the  following matrix operations.
  \begin{enumerate}
     \item[\textup{(i)}] $\circ: \Omega^{\mathbb{Z}\times \mathbb{Z}}\times \Omega^{\mathbb{Z}\times \mathbb{C}}\to
     \Omega^{\mathbb{Z}\times \mathbb{C}}$ by defining
     \begin{equation}\nonumber [A\circ X]_{m, x}:=\sum_{k}[A]_{m,k}[X]_{k, x};\end{equation}
 \item[\textup{(ii)}]  $\bullet_{(a,b)}: \Omega^{\mathbb{Z}\times \mathbb{C}}\times \Omega^{\mathbb{C}\times \mathbb{Z}}\to \Omega^{\mathbb{Z}\times \mathbb{Z}}$ by defining \begin{equation}\nonumber [X\bullet_{(a,b)} Y]_{m, n}:=\int_{a}^bW(x)[X]_{m,x}[Y]_{x,n}dx;\end{equation}
 \item[\textup{(iii)}]  $\bullet_{(q;a,b)} : \Omega^{\mathbb{Z}\times \mathbb{C}}\times \Omega^{\mathbb{C}\times \mathbb{Z}}\to \Omega^{\mathbb{Z}\times \mathbb{Z}}$ by
 defining \begin{equation}\nonumber [X\bullet_{(q;a,b)} Y]_{m, n}:=\int_{a}^bW(x)[X]_{m,x}[Y]_{x,n}d_qx;\end{equation}
 \item[\textup{(iv)}]  the transpose $X^T$ of $X$ by
$[X^T]_{x,m}:=[X]_{m,x}.$
 \end{enumerate}
In the above,   $W(x)$ is the general weight function, and the $q$-integral is defined by {\rm (cf.\cite[(1.11.2)-(1.11.3)]{10})}
\begin{align}\nonumber
\int_{b}^{a}f(t)d_qt&:=\int_{0}^{a}f(t)d_qt-\int_{0}^{b}f(t)d_qt, \\
\int_{0}^{a}f(t)d_qt&:=a(1-q)\sum_{k=0}^\infty f(aq^{k})q^k.\nonumber
\end{align}
\end{dy}
 A full and rigorous study on  Definition \ref{newmatrix}  will be given in our forthcoming paper. Here, we only point  out that
 the multiplications  $\circ$, $\bullet_{(a,b)}$ and $\bullet_{(q;a,b)}$ obey  the associative law, namely,
 \begin{xinzhi} Under the assumption of  Definition \ref{newmatrix},  there hold
 \begin{align}\nonumber
 A\circ (X\bullet_{(a,b)} Y)=(A\circ X)\bullet_{(a,b)} Y,~~
 A\circ (X\bullet_{(q;a,b)} Y)=(A\circ X)\bullet_{(q;a,b)} Y.
 \end{align}
 \end{xinzhi}
 With the help of these new matrix operations,  we are ready to consider two kinds of continuous $q$-orthogonal polynomials, as further illustrations of our idea.

\subsection{Dual form of the orthogonality relation of the big $q$-Jacobi
polynomials} It is well-known that  the big $q$-Jacobi polynomials (cf. \cite[(3.5.1)]{schem} or
\cite[(7.3.10)]{10}) are defined by
\begin{align}
P_n(x;a,b,c;q):=\,_3\phi_2\left[%
\begin{array}{ccc}
  q^{-n}, & abq^{n+1},& x\\
   & aq,&cq\\
\end{array};q,q
\right].\label{bjacobi}
\end{align}
From \cite[(3.5.2)]{schem},  it is clear
\begin{yl}
The big $q$-Jacobi polynomials
$P_n(x;a,b,c;q)$ satisfy the orthogonality relation
\begin{align}
\int_{cq}^{aq}\frac{\poq{x/a,x/c}{\infty}}{\poq{x,bx/c}{\infty}}P_n(x;a,b,c;q)P_m(x;a,b,c;q)d_qx
=h_n(a,b,c;q)\delta_{n,m},\label{731orthognew}
\end{align}
where
\begin{align}&
h_n(a,b,c;q)=M\frac{1-abq}{1-abq^{2n+1}}\frac{\poq{q,bq,abq/c}{n}}{
\poq{aq,cq,abq}{n}}(-acq^2)^{n}q^{n(n-1)/2}\label{H}\\
\qquad{and}&~~M=\frac{aq(1-q)\poq{q,c/a,aq/c,abq^{2}}{\infty}}{\poq{aq,bq,cq,abq/c}{\infty}}.\label{M}
\end{align}
\end{yl}
A proof of
(\ref{731orthognew})  can be found in  \cite{10} with the $q$-Gauss sum \cite[(II.8)]{10} involved. Instead here,  in light of inverse relations, we can show
\begin{dl} The dual form  of the orthogonality relation (\ref{731orthognew}) is  the following
transformation  of $\,_3\phi_2$ series
 \begin{align}
&\,_3\phi_2\left[%
\begin{array}{cccccccc}
   a, & b, &c\\
  & d, & e \\
\end{array};q,de/abc
\right]\nonumber\\
&=\frac{\poq{e/a,e/b,e/c}{\infty}}{\poq{e,e/d,de/abc}{\infty}}\,_3\phi_2\left[%
\begin{array}{cccccccc}
  d/a, &d/b,& d/c \\
  & d,&dq/e \\
\end{array};q,q
\right]\label{dl1166}\\ &+
\frac{\poq{d/a,d/b,d/c}{\infty}}{\poq{d,d/e,de/abc}{\infty}}\,_3\phi_2\left[%
\begin{array}{cccccccc}
 e/a,& e/b,&e/c \\
  & e,&eq/d \\
\end{array};q,q
\right],\nonumber
\end{align} where all parameters subject to convergent conditions of the sums.
    \end{dl}
\pf As indicated above, in order to express  (\ref{731orthognew}) in terms of matrix operations, we
 now  define the $q$-integral for  any matrix $\left(a_{i,j}(t)\right)$ to be
\begin{align}
\int_{0}^{a}\left(a_{i,j}(t)\right)d_qt:=\left(\int_{0}^{a}a_{i,j}(t)d_qt\right).\label{matrixintegral}
\end{align}
Besides, we need yet to introduce three matrices $A, X \in \Omega^{\mathbb{N}\times \mathbb{C}}$, $B\in \Omega^{\mathbb{C}\times \mathbb{C}}$ such that
\begin{align}\nonumber
 [A]_{n,t}=P_n(t;a,b,c;q), [X]_{n,t}=\frac{\poq{t,abq^2}{n}}{\poq{aq,cq}{n}} ,
\end{align}
and
\begin{align}\nonumber
[B]_{t,s}=\frac{\poq{t/a,t/c}{\infty}}{\poq{t,bt/c}{\infty}}\delta_{t,s},
\end{align}
where the (generalized) Kronecker delta $\delta_{t,s}=1$ for $t=s$ and $0$ otherwise.
As such, (\ref{bjacobi}) is now restated equivalently as
\begin{align}
 A=\textsf{N}(abq)\circ
  X.\label{sss}
 \end{align}
Next,  taking  \eqref{matrixintegral}  and the inverse relation $(\textsf{N}(a), \textsf{N}^{-1}(a))$ into account, we now  restate the
orthogonality relation  \eqref{731orthognew}  of  the big $q$-Jacobi polynomials
  as (in this case,  $W(x)=1$)
\begin{align}
 A\bullet_{(q;cq,aq)} B\bullet_{(q;cq,aq)} A^T=\left(h_n(a,b,c;q)\delta_{n,m}\right),
\end{align}
which, after  (\ref{sss}) inserted,  becomes
\begin{align}\nonumber
\textsf{N}(abq)\circ\left\{X\bullet_{(q;cq,aq)} B \bullet_{(q;cq,aq)} X^T\right\}\circ\textsf{N}(abq)^T=\left(h_n(a,b,c;q)\delta_{n,m}\right).
\end{align}
By inverting $\textsf{N}(abq)$ and $\textsf{N}(abq)^T$, we obtain
\begin{align}
X\bullet_{(q;cq,aq)} B\bullet_{(q;cq,aq)} X^T=\textsf{N}^{-1}(abq)\left(h_n(a,b,c;q)\delta_{n,m}\right)\textsf{N}^{-T}(abq).\label{tttt}
\end{align}
Next, substituting    \eqref{H} for $h_n(a,b,c;q)$ and equating the $(n,m)$-th entries on both sides of the  identity obtained, then we find immediately
\begin{align}
&\nonumber\int_{cq}^{aq}\frac{\poq{t/a,t/c}{\infty}}{\poq{t,bt/c}{\infty}}\frac{\poq{t,abq^2}{n}}{\poq{aq,cq}{n}}
\frac{\poq{t,abq^2}{m}}{\poq{aq,cq}{m}}d_qt\\
&\nonumber=M\sum_{k=0}^\infty
\frac{\poq{q^{-n}}{k}}{\poq{abq^{2+n}}{k}}\frac{\poq{abq}{k}}{\poq{q}{k}}\frac{1-abq^{2k+1}}{1-abq}q^{kn}
\frac{\poq{abq}{k}}{\poq{q}{k}}\frac{\poq{q^{-m}}{k}}{\poq{abq^{2+m}}{k}}\\
&\nonumber\frac{1-abq^{2k+1}}{1-abq}q^{km}\frac{1-abq}{1-abq^{2k+1}}\frac{\poq{q,bq,abq/c}{k}}{
\poq{aq,cq,abq}{k}}(-acq^2)^{k}q^{k(k-1)/2},
\end{align}
where $M$ is the same as in (\ref{M}).
After some rearrangement and simplification,  it reduces to
\begin{align}
&\lim_{d\to 0}\,_{8}W_{7}\left(abq;bq,
abq/c,acq^2/d,x,y;q,\frac{d}{xy}\right)\bigg|_{x=q^{-m},y=q^{-n}}
\nonumber\\
&=\frac{\poq{aq,bq,cq,abq/c}{\infty}}{aq(1-q)\poq{q,c/a,aq/c,abq^{2}}{\infty}}\frac{\poq{abq^2}{n}\poq{abq^2}{m}}{\poq{aq,cq}{n}\poq{aq,cq}{m}}
\nonumber\\
&\qquad\times\int_{cq}^{aq}\frac{\poq{t}{n}
\poq{t}{m}\poq{t/a,t/c}{\infty}}{\poq{t,bt/c}{\infty}}d_qt.\label{iddds}
\end{align}
In this case,  we can apply Watson's $q$-Whipple transformation (III.17) of  \cite{10}
to the $\,_8\phi_7$ series on the left side of (\ref{iddds}). The result is
\begin{align}
&\nonumber\,_{8}W_{7}\left(abq;bq,
abq/c,acq^2/d,x,y;q,\frac{d}{xy}\right)\bigg|_{x=q^{-m},y=q^{-n}}
\\
&\nonumber=\frac{(abq^2,bdq^{n}/c;q)_{m}}
{(bd/c,abq^{n+2};q)_{m}}\,_4\phi_3\left[%
\begin{array}{cccccccc}
  c/b, & acq^2/d, & q^{-n}, & q^{-m} \\
  & aq, & cq,& cq^{1-n-m}/bd \\
\end{array};q,q
\right].\nonumber
\end{align}
We therefore  get
\begin{align} &\int_{cq}^{aq}\frac{\poq{t}{n}
\poq{t}{m}\poq{t/a,t/c}{\infty}}{\poq{t,bt/c}{\infty}}d_qt\nonumber\\
&=
\frac{aq(1-q)\poq{q,c/a,aq/c}{\infty}\poq{abq^{2+m+n}}{\infty}\poq{aq,cq}{m}}{\poq{aq^{1+n},bq,cq^{1+n},abq/c}{\infty}}\label{111}\\
&\times\,_3\phi_2\left[%
\begin{array}{cccccccc}
  c/b,  & q^{-n}, & q^{-m} \\
  & aq, & cq \\
\end{array};q,abq^{2+n+m}
\right].\nonumber
\end{align}
It is easy to check that (\ref{111}) is equivalent to the following product of  three analytic functions in variables
$y,z$
\begin{align}F(y,z)=G(y,z)H(y,z)\label{finally}\end{align}
at the point $(y,z)=(q^n,q^m), m,n\geq 0$,  where  $F(y,z), G(y,z)$ and $H(y,z)$ are defined, respectively, by
\begin{align}
\nonumber F(y,z)&=\int_{cq}^{aq}\frac{\poq{t,t/a,t/c}{\infty}}{\poq{ty,tz,bt/c}{\infty}}d_qt,\\
\nonumber G(y,z)&=
\frac{aq(1-q)\poq{q,c/a,aq/c,abyzq^{2}, aq,cq}{\infty}}{\poq{aqy,aqz,cqy,cqz,bq,abq/c}{\infty}},\\
\nonumber H(y,z)&=\,_3\phi_2\left[%
\begin{array}{cccccccc}
  c/b,  & 1/y, & 1/z \\
  & aq, & cq \\
\end{array};q,abyzq^{2}
\right].
\end{align}
By analytic continuation,  we find that (\ref{finally}) holds as a
nonterminating version of (\ref{111}) on a disk around $(0,0)$. Reformulate
$F(y,z)$ on the left side of (\ref{finally}) in notation of  hypergeometric series, divide both sides
of the resulting identity  by $G(y,z)$, and finally make the replacement $(y,z)\to (1/y,1/z)$. Then  we obtain
\begin{align}
&\,_3\phi_2\left[%
\begin{array}{cccccccc}
   y, & z, &c/b\\
  & aq, & cq \\
\end{array};q,abq^{2}/yz
\right]\nonumber\\
&=\frac{\poq{cq/y,cq/z,bq}{\infty}}{\poq{cq,c/a,abq^{2}/yz}{\infty}}\,_3\phi_2\left[%
\begin{array}{cccccccc}
  aq/y, &aq/z,& abq/c \\
  & aq,&aq/c \\
\end{array};q,q
\right]\label{dl1122}\\ &+
\frac{\poq{aq/y,aq/z,abq/c}{\infty}}{\poq{aq,a/c,abq^{2}/yz}{\infty}}\,_3\phi_2\left[%
\begin{array}{cccccccc}
  cq/y,& cq/z,&bq \\
  & cq,&cq/a \\
\end{array};q,q
\right].\nonumber
\end{align}
 As a last step,  by making the simultaneous replacement of parameters
\begin{equation}\nonumber(a,b,c,y,z,q)\to (d/q,e/aq,e/q,b,c,q)\end{equation}
in (\ref{dl1122}), (\ref{dl1166}) follows at once. The theorem is now proved.\qed

More interestingly,  by putting   $b=c$ in (\ref{dl1122}),
we recover the nonterminating form of the $q$-Chu-Vandermonde formula.
\begin{tl}{\rm (cf.
\cite[(II.23)]{10})} For $a/c,y/c,z/c\not\in\{q^{-k}: 1\leq  k\in{\mathbb{N}}\}$,  it holds
\begin{align}
\,_2\phi_1\left[%
\begin{array}{cccccccc}
 aq/y,& aq/z \\
  & aq/c \\
\end{array};q,q
\right]
&+\frac{\poq{c/a,aq/y,aq/z}{\infty}}{\poq{a/c,cq/y,cq/z}{\infty}}\,_2\phi_1\left[%
\begin{array}{cccccccc}
  cq/y,& cq/z \\
  & cq/a \\
\end{array};q,q
\right]\nonumber\\
&=\frac{\poq{c/a,acq^{2}/yz}{\infty}}{\poq{cq/y,cq/z}{\infty}},\end{align}
 provided that the sums are convergent.
\end{tl}
Another two  special cases that $z=cq^{1+n}$ and $(y,z)=(aq^{1+m},cq^{1+n}),
m,n\geq 0,$ of (\ref{dl1122}) are worth mentioning.
\begin{tl}For $a,b,c\not\in\{q^{-k}: 1\leq k\in{\mathbb{N}}\}$ and $m,n\geq 0$,  there hold
\begin{align}
(\textup{i})\quad &\,_3\phi_2\left[%
\begin{array}{cccccccc}
  c/b,  & y, & cq^{1+n} \\
  & aq, & cq \\
\end{array};q,abq^{1-n}/yc
\right]\nonumber\\
&=\frac{\poq{aq/y,abq/c}{\infty}\poq{cq/a}{n}}{\poq{abq/yc,
aq}{\infty}\poq{yc/ab}{n}}\left(\frac{y}{bq}\right)^n\,_3\phi_2\left[%
\begin{array}{cccccccc}
 q^{-n}, & cq/y,&bq \\
  & cq,&cq/a \\
\end{array};q,q
\right];\nonumber\\
(\textup{ii})\quad&\,_3\phi_2\left[%
\begin{array}{cccccccc}
  c/b,  & aq^{1+m}, & cq^{1+n} \\
  & aq, & cq \\
\end{array};q,bq^{-m-n}/c
\right]= 0.
\end{align}
\end{tl}
\subsection{Dual form of the orthogonality relation of the  Askey-Wilson polynomials}
In what follows,  we will consider the orthogonality relation of the famous Askey-Wilson polynomials. As a member on the top of the Askey-scheme \cite{schem}, the Askey-Wilson
polynomials in a variable $x=\cos\theta$ is defined by \cite[(3.1.1)]{schem}
\begin{align}
W_n(x;a,b,c,d|q):=\frac{\poq{ab,ac,ad}{n}}{a^n}\,_4\phi_3\left[%
\begin{array}{cccc}
  q^{-n}, & abcdq^{n-1},&ae^{i\theta},& ae^{-i\theta}\\
   & ab,&ac,&ad\\
\end{array};q,q
\right].\label{askey}
\end{align}
It is well known  that
\begin{yl}{\rm (cf. \cite[(3.1.2)]{schem})}
 The Askey-Wilson polynomials satisfy the orthogonality relation
\begin{align}
\frac{1}{2\pi}\int_{-1}^1\frac{w(x)}{\sqrt{1-x^2}}W_n(x;a,b,c,d|q)W_m(x;a,b,c,d|q)dx
=\kappa_n\delta_{n,m},\label{738orthog2}
\end{align}
where
\begin{align}
w(x)&=\frac{h(x,1)h(x,-1)h(x,q^{1/2})h(x,-q^{1/2})}{h(x,a)h(x,b)h(x,c)h(x,d)},\label{added1}\\
h(x,y)&=\poq{ye^{i\theta},
ye^{-i\theta}}{\infty},\label{added2}\\
\kappa_n&=\frac{\poq{abcdq^{n-1}}{n}\poq{abcdq^{2n}}{\infty}}
{\poq{q^{n+1},abq^n,acq^n,adq^n,bcq^n,bdq^n,cdq^n}{\infty}}.\nonumber
\end{align}
\end{yl}
Now we are able to present  the dual form of (\ref{738orthog2}) as the following theorem, whose derivation may be regarded as a new proof of (\ref{738orthog2}) via the integral representation of VWP $\,_8\phi_7$ series \cite[\S 6.3]{10}.

\begin{dl}\label{dl1144} For $\max\{|a|,|b|,|c|,|d|\}<1$ and $a,b,c,d\in \mathbb{C}$, the dual form  of the orthogonality relation (\ref{738orthog2}) is  the   Askey-Wilson $q$-beta integral  \begin{align}
&\frac{1}{2\pi}\int_{-1}^1\frac{w(x)}{\sqrt{1-x^2}}\poq{ae^{i\theta},ae^{-i\theta}}{n}\poq{ae^{i\theta},ae^{-i\theta}}{m}
dx\label{5t}\\
&=\mathcal{K}\frac{\poq{ab,ac,ad}{m}\poq{ab,ac,ad}{n}}{\poq{abcd}{m}\poq{abcd}{n}}\,_{8}W_{7}\left(abcd/q;bc,bd,cd,q^{-m},q^{-n};q,a^2q^{m+n}\right),\nonumber\end{align}
where
\begin{equation}\nonumber\mathcal{K}=\frac{\poq{abcd}{\infty}}{\poq{q,ab,ac,ad,bc,bd,cd}{\infty}}.\end{equation}
        \end{dl}
\pf It suffices to define three matrices $A,D,$ and $X$, respectively, by
\begin{align}\nonumber
[A]_{n,x}=W_n(x;a,b,c,d|q),  [D]_{n,k}=\frac{\poq{ab,ac,ad}{n}}{a^{n}}\delta_{n,k},
\end{align}
and
\begin{equation}\nonumber[X]_{n,x}=\frac{\poq{abcd,ae^{i\theta},ae^{-i\theta}}{n}}{\poq{ab,ac,ad}{n}}.\end{equation}
Further,  we  choose the weight function
\begin{equation}\nonumber W(x)=\frac{1}{2\pi}\frac{w(x)}{\sqrt{1-x^2}},\end{equation}
where   $x\in [-1,1]$.  In light of the matrix
operation $\circ$,
   it is easily found that (\ref{askey})
 is just equivalent to
 \begin{equation}\nonumber A=D\times\textsf{N}(abcd/q)\circ X,\end{equation}
 where  ``$\times$" denotes the usual  matrix multiplication. Thus the orthogonality relation (\ref{738orthog2}),  by use of the matrix operation $\bullet_{(-1,1)} $,
 can still be restated equivalently as
\begin{equation}\nonumber A\bullet_{(-1,1)} A^T=\left(\kappa_n\delta_{n,m}\right),\end{equation} viz., \begin{align}D\times\textsf{N}(abcd/q)\circ X\bullet_{(-1,1)} X^T\circ \textsf{N}(abcd/q)^T\times D^T=\left(\kappa_n\delta_{n,m}\right),
\end{align}
from which we may deduce \begin{align}
X\bullet_{(-1,1)} X^T=\textsf{N}^{-1}(abcd/q)\times\left(\frac{a^{2n}\kappa_n}{\poq{ab,ac,ad}{n}^2}
\delta_{n,m}\right)\times
\textsf{N}(abcd/q)^{-T}.\label{ttt2}
\end{align}
Thus, by equating the $(n,m)$-th entries on both sides  of
(\ref{ttt2}), we  obtain  \eqref{5t}. \qed

It is worth pointing out that (\ref{5t}) is in agreement with (6.3.8) of \cite{10}  under the simultaneous
 replacement of parameters
\begin{equation}\nonumber(a,b,c,d,f,g)\to (b,c,d,aq^m,aq^n,a).\end{equation} From a quick glance over \cite{nassrahman}, we know that   Nassrallah and Rahman  proved (6.3.8) of \cite{10} via the direct use of $q$-Askey-Wilson beta integral without  appeal to  the orthogonality relation of the Askey-Wilson polynomials $W_n(x;a,b,c,d|q)$.  It is also of interest to see that the limitation $m,n\to\infty$ of Theorem \ref{dl1144}  yields

\begin{tl}\label{dl114455} Let $w(x)$ and $h(x,a)$ be given by \eqref{added1} and \eqref{added2}, respectively. Then
 \begin{align}
\int_{-1}^1\frac{w(x)}{\sqrt{1-x^2}}h^2(x,a)dx&=\frac{2\pi\poq{ab,ac,ad}{\infty}}{\poq{q,bc,bd,cd,abcd}{\infty}}\label{5s}
\\&\times\sum_{k=0}^{\infty}\frac{1-abcdq^{2k-1}}{1-abcd/q}\frac{\poq{abcd/q,bc,bd,cd}{k}}{\poq{q,ab,ac,ad}{k}}
a^{2k}q^{k(k-1)}.\nonumber
\end{align}
        \end{tl}

 \section{Concluding remarks}
So far, we have applied the inverse relation $(\textsf{N}(a), \textsf{N}^{-1}(a))$ to find the dual forms of the orthogonality relations for some discrete and continuous $q$-orthogonal polynomials.  Apart from  $\textsf{N}(a)$, there are many  inverse relations can be used. For instance, let $\textsf{M}(a)$ and $\textsf{K}(a,c)$ be two ILT matrices   with the  entries given by
\begin{align}
\nonumber[\textsf{M}(a)]_{n,k}&=\frac{\poq{q^{-n},1/a}{k}}{\poq{q,q^{1-n}/a}{k}}q^k\qquad\mbox{and}\\
\nonumber[\textsf{K}(a,c)]_{n,k}&=\frac{\poq{aq^n, q^{-n},q\sqrt{c},-q\sqrt{c},c/a,c}{k}}{\poq{cq^{1-n}/a,
cq^{n+1},q,\sqrt{c},-\sqrt{c},aq}{k}}q^k,
\end{align}
respectively.
 Then it is easy to verify by Theorem \ref{math4} that
 \begin{align}\nonumber
 \textsf{M}^{-1}(a)=\textsf{M}(1/a),~~ \textsf{K}^{-1}(a,c)=\textsf{K}(c,a).
 \end{align}
 We remark that $\textsf{K}^{-1}(a,c)$  is consistent with the
Krattenthaler inversion formula (cf. \cite[(1.5)(1)/(1.5)(2)]{kra}). Without any doubt, the potential applications of these two inverse relations similar  to those of $\textsf{N}(a)$  deserve further investigation. Generally speaking,   inverse relations are helpful  to reveal possibly association of orthogonality relations like (\ref{lagg})/(\ref{laggorth}) of  discrete and continuous $q$-orthogonal polynomials with summation and transformation of hypergeometric series, even more in solving  the problem of connection coefficients between $q$-orthogonal polynomials.

\bibliographystyle{amsplain}

\end{document}